\documentclass[reqno, 12pt]{amsart}
\def\C{{\mathbb C}}

\def\Q{{\mathbb Q}}
\def\Qp{{\mathbb Q}_p}

\def\Zp{{\mathbb Z}_p}

\def\F{{\mathbb F}}
\def\Qp {{{\mathbb Q}_p}}
\def\Fq{{{\mathbb F}_q}}

\newtheorem{cor}{Corollary}
\newtheorem{prop}{Proposition}

\newtheorem{theorem}{Theorem}

\theoremstyle{definition}
\newtheorem{defn}{Definition}

\theoremstyle{remark}
\newtheorem{rem}{Remark}

\begin{document}
\title[The Ongoing Binomial Revolution]
{The Ongoing Binomial Revolution}
\author{David Goss}
\thanks{}
\address{Department of Mathematics\\The Ohio State University\\231
W.\ $18^{\rm th}$ Ave.\\Columbus, Ohio 43210}

\email{goss@math.ohio-state.edu}
\date{Spring, 2011}

\begin{abstract}
The Binomial Theorem has long been essential in 
mathematics. In one form or another it was known to the
ancients and, in the hands of Leibniz, Newton, Euler, Galois, and others, it became
an essential tool in both algebra and analysis. Indeed, Newton early on developed certain binomial
series (see Section \ref{newton}) which played a role in his subsequent work
on the calculus. From the work of Leibniz, Galois, Frobenius, and many others, 
we know of its
essential role in algebra. In this paper we rapidly trace the history  of the Binomial Theorem,
binomial series, and binomial coefficients, with emphasis on their decisive role in
function field arithmetic. We also explain conversely
how function field arithmetic is now leading to
new results in the binomial theory via insights into characteristic $p$ $L$-series.
\end{abstract}

\maketitle

\section{Introduction}\label{intro}
The Binomial Theorem has played a crucial role in the development of mathematics,
algebraic or analytic, pure or applied. It was very important in the development
of the calculus, in a variety of ways, and has certainly been as important in the 
development of number theory. It plays a dominant
role  in function field arithmetic. In fact, it 
almost appears as if function field arithmetic ({\em and} a large chunk of arithmetic
in general) is but a commentary on this amazing result. In turn, function field
arithmetic has recently returned the favor by shedding new light on the Binomial
Theorem. It is our purpose here to recall the history of the Binomial Theorem, with
an eye on applications in characteristic $p$, and finish by discussing these
new results.

We obviously make no claims here to being encyclopedic. Indeed, to thoroughly
cover the Binomial Theorem would take many volumes. Rather, we have chosen to walk
a quick and fine line through the many relevant results.

This paper constitutes a serious reworking of my April 2010 lectures at the 
{\it Centre de Recerca Matem\`atica} in Barcelona. It is my great pleasure to thank
the organizers of the workshop and, in particular, Francesc Bars.

\section{Early History}\label{early}
According to our current understanding, the Binomial Theorem can be
traced to the $4$-th century B.C. and Euclid where one finds the formula for $(a+b)^2$.
In the $3$-rd century B.C. the Indian mathematician Pingala presented what is
now known as ``Pascal's triangle'' giving binomial coefficients in a triangle. 
Much later, in the 10-th century A.D., the Indian mathematician Halayudha
and the Persian mathematician al-Karaji derived similar results as did the
13-th century Chinese mathematician Yang Hui. It is remarkable that al-Karaji appears
to have used mathematical induction in his studies.

Indeed, binomial coefficients, appearing in Pascal's triangle, seem to have been
widely known in antiquity. Besides the mathematicians mentioned above, Omar
Khayyam (in the 11-th century), Tartaglia, Cardano, Vi\'ete, Michael Stifel (in the
16-th century), and William Oughtred, John Wallis, Henry Briggs, and Father 
Marin Mersenne (in the 17-th century) knew of these numbers.
In the 17-th century, Blaise Pascal gave the binomial coefficients their now commonly
used form: for a nonnegative integer $n$ one sets
\begin{equation}\label{pascal}
{n \choose k}:=\frac{n!}{k!(n-k)!}=\frac{n(n-1)(n-2)\cdots(n-k+1)}
{k(k-1)(k-2)\cdots1}\,.
\end{equation}
With this definition we have the very famous, and equally ubiquitous, 
 {\it Binomial Theorem}: 
\begin{equation}\label{binomial2} (x+y)^n=\sum_k {n \choose k}x^k y^{n-k}\,.
\end{equation}
And of course, we also deduce the first miracle giving the integrality of the
binomial coefficients $n \choose k$.

Replacing $n$ by a variable $s$ in Equation \ref{pascal} gives the
{\it binomial polynomials} 
\begin{equation}\label{binomial1}
{s\choose k}:= \frac{s(s-1)(s-2)\cdots(s-k+1)}{k(k-1)(k-2)\cdots1}=
\frac{s(s-1)(s-2)\cdots(s-k+1)}{k!}\,.
\end{equation}

\section{Newton, Euler, Abel, and Gauss}\label{newton}
We now come to Sir Isaac Newton and his contribution to the Binomial Theorem.
His contributions evidently were discovered in the year 1665 (while sojourning
in Woolsthorpe, England to avoid an outbreak of the plague) and 
discussed in a letter to Oldenburg in 1676. Newton was highly influenced by
work of John Wallis who was able to calculate the area under the curves
$(1-x^2)^n$, for $n$ a nonnegative integer. Newton then considered 
fractional exponents $s$ instead of $n$. He realized that one could find the successive
coefficients $c_k$  of $(-x^2)^k$, in the expansion of $(1-x^2)^s$, by multiplying the
previous coefficient by $\frac{s-k+1}{k}$ exactly as in the integral case. In particular,
Newton formally computed the Maclaurin
series for $(1-x^2)^{1/2}$, $(1-x^2)^{3/2}$ and $(1-x^2)^{1/3}$.

(One can read about this in the paper \cite{Co1} where the author believes that
Newton's contributions to the Binomial Theorem were relatively minor and that
the credit for discussing fractional powers should go to James Gregory -- who in
1670 wrote down the series for $b\left(1+\frac{d}{b}\right)^{a/c}$. This is a distinctly
minority viewpoint.)

In any case, Newton's work on the Binomial Theorem played a role in his subsequent
work on calculus. However, Newton did not consider issues of convergence. 
This was discussed by Euler, Abel, and Gauss. Gauss gave the first satisfactory
proof of convergence of such series in 1812.  Later Abel gave a treatment that would
work for general complex numbers. The theorem on binomial series can now
be stated.
\begin{theorem}\label{binomial3}
Let $s\in \C$. Then the series $\sum_{k=0}^\infty {s \choose k} x^k$ converges to
$(1+x)^s$ for all complex $x$ with $\vert x\vert<1$.\end{theorem}

\begin{rem}\label{binomial4} Let $s=n$ be any integer, positive or negative. 
Then for all complex
$x$ and $y$ with $\vert x/y\vert<1$ one readily deduces from Theorem
\ref{binomial3}  a convergent expansion
\begin{equation}\label{binomial5}
(x+y)^n=\sum_{k=0}^\infty {n \choose k}x^ky^{n-k}\,.\end{equation}\end{rem}

It is worth noting that Gauss' work on the convergence of the binomial series
marks the first time convergence involving {\it any} infinite series was satisfactorily
treated!

Now let $f(x)$ be a polynomial with coefficients in an extension of $\Q$.
The degree of $x \choose k$ as a polynomial in $x$ is $k$. As such one
can always expand $f(x)$ as a linear combination of $x \choose k$.  Such
an expansion is called the {\it Newton series} and can be traced back to
his {\it Principia Mathematica} (1687). The coefficients of such an expansion are
given as follows.
\begin{defn}\label{binomial6} Set $(\Delta f)(x):=f(x+1)-f(x)$. \end{defn}
\begin{prop}\label{binomial7} We have 
\begin{equation}\label{binomial8}
f(x)=\sum_k(\Delta^k f)(0) {x \choose k}\,.\end{equation}\end{prop}

\section{The $p$-th power mapping}\label{pthpower} 
Let $p$ be a prime number and let $\F_p$ be the field with $p$-elements.
The following elementary theorem is then 
absolutely fundamental for number theory and arithmetic
geometry. Indeed its importance cannot be overstated.
\begin{theorem}\label{pthpower1}
Let $R$ be any $\F_p$-algebra. Then the mapping $x\mapsto x^p$ is
a homomorphism from $R$ to itself.\end{theorem}
\noindent
As is universally known, the proof amounts to expanding by the Binomial
Theorem and noting that for $0<i<p$, one
has ${p \choose i}\equiv 0\pmod{p}$ as the denominator of Equation \ref{pascal}
is prime to $p$.

According to Leonard Dickson's history (Chapter III of \cite{Di1}), the first person
to establish (a form of) Theorem \ref{pthpower1} was 
Gottfried Leibniz on September 20, 1680. One can then rapidly deduce a proof of
Fermat's Little Theorem (i.e., $a^p\equiv a \pmod{p}$ for all integers $a$ and primes
$p$). Around 1830 Galois used interates of the $p$-th power mapping to construct
general finite fields.

It was 216 years after Leibniz (1896) that the
equally essential {\it Frobenius automorphism} (or {\it Frobenius substitution}) in
the Galois theory of fields was born. Much of modern number theory and algebraic
geometry consists of computing invariants of the $p$-th power mapping/Frobenius map.

Drinfeld modules are subrings of the algebra (under composition!) of ``polynomials''
in the $p$-th power mapping; thus their very existence depends on the 
Binomial Theorem.

\begin{rem}\label{pthpower2} With regard to the proof of Theorem \ref{pthpower1},
it should also be noted that Kummer in 1852 established that the exact power
of a prime $p$ dividing ${n \choose k}$ is precisely the number of ``carries''involved
in adding $n-k$ and $k$ when they are expressed in their canonical $p$-adic
expansion.\end{rem}
\section{The Theorem of Lucas}\label{lucas}
Basic for us, and general arithmetic in finite characteristic, is the famous Theorem of
Lucas from 1878 \cite{Lu1}. Let $n$ and $k$ be two nonnegative integers and
$p$ a prime. Write $n$ and $k$ $p$-adically as $n=\sum_i n_ip^i$, $0\leq n_i<p$
and $k=\sum_i k_ip^i$, $0\leq k_i<p$.
\begin{theorem}\label{lucas1} {\bf (Lucas)}
We have 
\begin{equation}\label{lucas2}
{n\choose k}=\prod_i {n_i \choose k_i}\pmod{p}\,.\end{equation}
\end{theorem}
\begin{proof}
We have $(1+x)^n=(1+x)^{\sum n_ip^i}=\prod_i(1+x)^{n_ip^i}$. Modulo $p$, Theorem 
\ref{pthpower} implies that $(1+x)^n=\prod_i (1+x^{p^i})^{n_i}$. The result then
follows by expressing both sides by the Binomial Theorem and
the uniqueness of $p$-adic expansions.\end{proof}
\section{The Theorem of Mahler}\label{Mahler}
The binomial polynomials $s \choose k$ (given in Equation\ref{binomial1}) obviously
have coefficients in $\Q$ and thus also can be considered in the $p$-adic numbers
$\Qp$.
\begin{prop}\label{Mahler1}
The functions $s \choose k$, $k=0,1,\ldots$, map $\Zp$ to itself.
\end{prop}
\begin{proof} Indeed, $s \choose k$ takes the nonnegative integers to themselves.
As these are dense in $\Zp$, and $\Zp$ is closed, the result follows. \end{proof}
  
Let $y\in \Zp$, and formally set $f_y(x):=(1+x)^y$. By the above proposition, 
$f_y(x)\in \Zp[[x]]$. As such, we can consider $f_y(x)$ in {\em any} 
nonArchimedean field of {\em any} characteristic where it will converge on
the open unit disc.

Let $\{a_k\}$ be a collection of $p$-adic numbers approaching $0$ as $k\to 
\infty$ and put $g(s)=\sum_k a_k{s \choose k}$; it is easy to see that this series
converges to a continuous function from $\Zp$ to $\Qp$.  
Moreover, given a continuous function $f\colon \Zp\to \Qp$, the Newton series
(Equation \ref{binomial8}) certainly makes sense formally.

\begin{theorem}\label{Mahler2} {\bf (Mahler)} The Newton series of a
continuous function $f\colon \Zp\to \Qp$ uniformly converges to it.
\end{theorem}
\noindent
The proof can be found in \cite{Ma1} (1958). 
The Mahler expansion of a continuous
$p$-adic function is obviously unique. 

Mahler's Theorem can readily be 
extended to continuous functions of $\Zp$ into complete fields of characteristic
$p$.  One can also find analogs of it that work for functions on the maximal
compact subrings of arbitrary local fields. In characteristic $p$, an especially
important analog of the binomial polynomials was constructed by L. Carlitz
as a byproduct of his construction of the Carlitz module (see, e.g., \cite{Wa1}). 

Carlitz's construction can be readily described. Let $e_k(x):=\prod (x-\alpha)$
where $\alpha$ runs over elements of $\Fq[t]$, $q=p^{m_0}$, of degree $<k$. 
As these elements
form a finite dimensional $\Fq$-vector space, the functions $e_k(x)$ are readily
seen to be $\Fq$-linear. Set $D_k:=e_k(t^k)=\prod f$ where $f(t)$ runs through
the {\em monic} polynomials of degree $k$. Carlitz then establishes that
$e_k(g)/D_k$ is {\it integral} for $g\in \Fq[t]$. 
\begin{rem}\label{Mahler2.5}
The binomial coefficients $s \choose k$ appear in the power series
expansion of $(1+x)^s$. It is very important to note that the the polynomials
$e_k(x)/D_k$ appear in a completely similar fashion in terms of the 
expansion of the {\it Carlitz module} -- an $\Fq[t]$-analog of $\mathbb G_m$;
see, e.g., Corollary 3.5.3 of \cite{Go1}.
\end{rem}

\noindent
Now let $k$ be any nonnegative integer written
$q$-adically as $\sum_t k_tq^t$, $0\leq k_t<q$ for all $t$.
\begin{defn}\label{Mahler3}
We set
\begin{equation}\label{Mahler4}
G_k(x):=\prod_t \left ( \frac{e_t(x)}{D_t}\right)^{k_t}\,.
\end{equation}\end{defn}
The set $\{G_k(x)\}$ is then an excellent characteristic $p$ replacement for
$\{{s \choose k}\}$ in terms of analogs of Mahler's Theorem, etc, see \cite{Wa1}. In 2000
K.\ Conrad \cite{Con1} showed that Carlitz's use of digits in constructing analogs
of $s \choose k$ can be applied quite generally.

In a very important refinement of Mahler's result, in 1964 Y.\ Amice \cite{Am1}  gave
necessary and sufficient conditions on the Mahler coefficients guaranteeing 
that a function can be locally expanded in power series. In fact, Amice's results
work for arbitrary local fields and are also essential for the function field theory.
Indeed, as the function $(1+x)^y$, $y\in \Zp$, is clearly locally analytic, Amice's
results show that its expansion coefficients tend to $0$ very quickly, thus allowing
for general analytic continuation of $L$-series and partial $L$-series 
\cite{Go2}.

In 2009, S.\ Jeong \cite{Je1} established that the functions $u\mapsto u^y$, $y\in \Zp$
precisely comprise the group of {\em locally-analytic} endomorphisms of the
$1$-units in a local field of finite characteristic.

\section{Measure Theory}\label{measure}
Given a local field $K$ with maximal compact $R$, one is able to describe a theory
of integration for all continuous $K$-valued functions on $R$.  
A {\it measure on $R$ with values in $R$} is a finitely-additive, $R$-valued
function on the compact open subsets of $R$. Given a measure $\mu$ and
a continuous $K$-valued function $f$ on $R$, the Riemann sums for
$f$ (in terms of compact open subsets of $R$) are easily seen to
converge to an element of $K$ naturally denoted $\int_R f(z)\, d\mu(z)$.

Given two measures $\mu_1$ and $\mu_2$, we are able to form their
convolution $\mu_1\ast \mu_2$ in exactly the same fashion as in classical
analysis. In this way, the space of measures forms a commutative $K$-algebra.

In the case of 
$\Qp$ and $\Zp$ one is able to use Mahler's Theorem (Theorem \ref{Mahler2} above)
to express integrals of general continuous functions in terms of the
integrals of binomial coefficients. 

Now $(1+z)^{x+y}=(1+z)^x (1+z)^y$ giving an {\it addition formula} for the
binomial coefficients. Using this in the convolution
 allows one to establish that the convolution algebra of measures (the {\it Iwasawa
algebra}) is isomorphic to $\Zp[[X]]$.

In finite characteristic, we obtain a {\em dual} characterization of measures that is
still highly mysterious and {\em also} depends crucially on the Binomial Theorem. So let
$q$ be a power of a prime $p$ as above. Let $n$ be a nonnegative integer written
$q$-adically as $\sum n_kq^k$. Thus, in characteristic $p$, we deduce
\begin{equation}\label{measure1}
(x+y)^n=\prod_k(x+y)^{n_kq^k}=\prod_k (x^{q^k}+y^{q^k})^{n_k}\, \end{equation}

Now recall the definition of the functions $G_n(x)$ (Definition \ref{Mahler3} above)
via digit expansions. As the functions $e_j(x)$ are also {\it additive} we immediately
deduce from Equation \ref{measure1} the next result.
\begin{theorem}\label{measure2}
We have
\begin{equation}\label{measure3}
G_n(x+y)=\sum_{j=0}^n {n \choose j} G_j(x)G_{n-j}(x)\,.\end{equation}\end{theorem}
\noindent
In other words, the functions $\{G_n(x)\}$, {\em also} satisfy the Binomial Theorem!

Let $\mathfrak D_j$ be the hyperdifferential (= ``divided derivative'') operator given
by ${\mathfrak D}_jz^i:={i \choose j}z^{i-j}$.  Notice that ${\mathfrak D}_i
{\mathfrak D}_j={i+j \choose i} {\mathfrak D}_{i+j}$. Let $R\{\{{\mathfrak D}\}\}$ be
the algebra of formal power-series in the $\mathfrak D_i$ with the above multiplication
rule where $R$ is any commutative ring. Note further that this definition makes
sense for all $R$ precisely since $i \choose j$ is always integral.

Let $A=\Fq[t]$ and let $f\in A$ be irreducible; set $R:=A_f$, the completion of
$A$ at $(f)$. Using the Binomial Theorem for the Carlitz polynomials we have
the next result \cite{Go3}

\begin{theorem}\label{measure4}
The convolution algebra of $R$-valued measures on $R$ is isomorphic to
$R\{\{{\mathfrak D}\}\}$.\end{theorem}

\begin{rem}\label{measure5}  The history of Theorem \ref{measure4} is amusing. 
I had calculated the algebra of measures using the Binomial Theorem and then 
showed the
calculation to Greg Anderson who, rather quickly(!), recognized it as the
ring of hyperderivatives/divided power series.\end{rem}
\begin{rem}\label{measure6} One can ask why we represent the algebra of
measures as operators as opposed to divided power series. Let $\mu$ be
a measure on $R$ ($R$ as above) and let $f$ be a continuous function; one can
then obtain a new continuous function $\mu(f)$ by 
\begin{equation}\label{measure7}
\mu(f)(x):=\int_R f(x+y)\, d\mu(y)\,.\end{equation}
The operation of passing from the expansion of $f$ (in the Carlitz
polynomials) to the expansion of
$\mu(f)$ formally appears as if the differential operator attached to $\mu$ acted
on the expansion. This explains our choice.\end{rem}

\section{The group $S_{(p)}$ and binomial symmetries in finite 
characteristic}\label{symm}
Let $q=p^{m_0}$, $p$ prime, as above, and let $y\in \Zp$. Write $y$ $q$-adically 
as
\begin{equation}\label{qadic}
y=\sum_{k=0}^\infty y_k q^k
\end{equation}
where $0\leq y_k<q$ for all $k$. If $y$ is a nonnegative integer 
(so that the sum in Equation
\ref{qadic} is obviously finite), then we set $\ell_q(y)=\sum_k y_k\,.$

Let $\rho$ be a permutation of the set $\{0,1,2,\ldots\}$.

\begin{defn}\label{pi(n)}
We define $\rho_\ast (y)$, $y\in \Zp$, by
\begin{equation}\label{pi(n)2} 
\rho_\ast(y):=\sum_{i=0}^\infty y_k q^{\rho(i)}\,.
\end{equation}\end{defn} 
Clearly $y\mapsto \rho_\ast(y)$ is a bijection of $\Zp$. Let $S_{(q)}$ be the
group of bijections of $\Zp$ obtained this way.
Note that if $q_0$ and $q_1$ are powers of $p$, and $q_0\mid q_1$, then
$S_{(q_1)}$ is naturally realized as a subgroup of $S_{(q_0)}$. 

\begin{prop}\label{basicS}
Let $\rho_\ast(y)$ be defined as above.\\
{\rm 1}. The mapping $y\mapsto \rho_\ast(y)$ is a homeomorphism of $\Zp$.\\
{\rm 2.} (``Semi-additivity'') Let $x,y,z$ be three $p$-adic integers
with $z=x+y$ and where there is no carry over of $q$-adic digits. Then
$\rho_\ast(z)=\rho_\ast(x)+\rho_\ast(y)$.\\
{\rm 3.} The mapping $\rho_\ast(y)$ stabilizes both the nonnegative and
nonpositive integers.\\
{\rm 4.} Let $n$ be a nonnegative integer. Then $\ell_q(n)=\ell_q(\rho_\ast(n))$.\\
{\rm 5}. Let $n$ be an integer. Then $n\equiv \rho_\ast(n) \pmod{q-1}$.\\
\end{prop}
\noindent
For the proof, see \cite{Go4}.
\begin{prop}\label{symm1}
Let $\sigma\in S_{(p)}$, $y\in \Zp$, and $k$ a nonnegative integer. Then
we have
\begin{equation}\label{symm2}
{y \choose k}\equiv  {\sigma y\choose \sigma k}\pmod{p}\,.\
\end{equation}
\end{prop}
\begin{proof} This follows immediately from the Theorem of Lucas
(Theorem \ref{lucas1}).\end{proof}
\begin{cor}\label{symm3} Modulo $p$, we have ${\sigma y \choose k}={y \choose
\sigma^{-1} k}$.\end{cor}
\begin{cor}\label{symm4} We have $p\mid {y\choose k} \Longleftrightarrow p\mid
{\sigma y \choose \sigma k}$.\end{cor}

\begin{prop}\label{symm5}
Let $i$ and $j$ be two nonnegative integers. Let $\sigma\in S_{(p)}$. Then
\begin{equation}\label{symm6}
{{i+j}\choose i}\equiv {{\sigma i+\sigma j}\choose \sigma i}\pmod{p}\,.
\end{equation}
\end{prop}
\begin{proof}
The theorems of Lucas and Kummer show that if there is any carry over of $p$-adic digits
in the addition of $i$ and $j$, then ${{i+j}\choose i}$ is $0$ modulo
$p$. However, there is carry over of the $p$-adic digits in the sum
of $i$ and $j$ if and only if there is carry over in the sum of
$\sigma i$ and $\sigma j$; in this case both sums are $0$ modulo $p$.
If there is no carry over, then the result follows from Part 2 of Proposition
\ref{basicS} and Proposition \ref{symm1}.\end{proof} 
Let $R$ be as in the previous section.
\begin{cor}\label{symm7}
The mapping $\mathfrak D_i\mapsto \mathfrak D_{\sigma i}$ is an automorphism
of $R\{\{\mathfrak D\}\}$. \end{cor}

It is quite remarkable that the group $S_{(q)}$ very much {\em appears} to be
a symmetry group of characteristic $p$ $L$-series. Indeed, in examples, this
group preserves the orders of trivial zeroes as well as the denominators of
special zeta values (the ``Bernoulli-Carlitz'' elements). Moreover, given a nonnegative
integer $i$, one has the ``special polynomials''  of characteristic $p$ $L$-series
arising at $-i$. It is absolutely
remarkable, and highly nontrivial to show, that the degrees of these special polynomials
are {\em invariant} of the action of $S_{(q)}$ on $i$. Finally, the action of
$S_{(q)}$ even appears to extend to the {\em zeroes} themselves 
of these characteristic $p$ functions. See \cite{Go4} for all this.

\section{The Future}\label{future}
We have seen how the Binomial Theorem has impacted the development of both
algebra and analysis. In turn these developments have provided the foundations
for characteristic $p$ arithmetic. Furthermore, as in Section \ref{symm},
characteristic $p$ arithmetic
has contributed results relating to the Binomial Theorem of both
an algebraic (automorphisms of $\Zp$ and binomial coefficients) and analytic
(automorphisms of algebras of divided derivatives) nature. Future
research should lead to a deeper understanding of these recent offshoots of the Binomial
Theorem as well as add many, as yet undiscovered, new ones,

\end{document}